\documentclass[12pt,leqno]{amsart}
\usepackage{amsmath}
\usepackage{amssymb}
\def\overset#1#2{{\mathrel{\mathop {{#2}_{}}\limits^{#1}}}}
\def\underset#1#2{{\mathrel{\mathop {{}_{} {#2}}\limits_{{#1}_{}}}}}
\def\upplim_#1{\underset{#1}{\overline\lim}\;}
\def\lowlim_#1{\underset{#1}{\underline\lim}\;}
\newtheorem{corollary}[equation]{Corollary}
\newtheorem{definition}[equation]{Definition}
\newtheorem{claim}[equation]{\indent{\it Claim}\rm }

\newtheorem{lem}[equation]{Lemma}
\newtheorem{lemma}[equation]{Lemma}
\newtheorem{proposition}[equation]{Proposition}

\newtheorem{theorem}[equation]{Theorem}

\newcommand{\C}{{\mathbf{C}}}

\renewcommand{\P}{{\mathbf{P}}}

\newcommand{\R}{{\mathbf{R}}}

\newcommand{\rank}{\mathrm{rank}}

\newcommand{\supp}{\mathrm{Supp}\,}

\newcommand{\Z}{\mathbf{Z}}

\numberwithin{equation}{section}

\title[Second main theorems for meromorphic mappings]{Second main theorems for meromorphic mappings and moving hyperplanes with truncated counting functions} 

\date { }
\author{Si Duc Quang}

\begin{document}

\begin{abstract}
In this article, we establish some new second main theorems for meromorphic mappings of $\C^m$ into $\P^n(\C)$ and moving hyperplanes with truncated counting functions. Our results are improvements of the previous second main theorems for moving hyperplanes with truncated (to level $n$) counting functions. 
\end{abstract}

\def\thefootnote{\empty}
\footnotetext{
2010 Mathematics Subject Classification:
Primary 32H30, 32A22; Secondary 30D35.\\
\hskip8pt Key words and phrases: Nevanlinna theory, second main theorem, meromorphic mapping, moving hyperplane.}

\maketitle

\section{Introduction}
The second main theorem for meromorphic mappings into projective spaces with moving hyperplanes was first given by W. Stoll, M. Ru \cite{RS} and M. Shirosaki in 1990's \cite{S1,S2}, where the counting functions are not truncated. In 2000, M. Ru \cite{MR} proved a second main theorem with trucated counting functions for nondegenerate mappings of $\C$ into $\P^n(\C)$ and moving hyperplanes. After that, this result was reproved for the case of several complex variables by Thai-Quang \cite{TQ05}. For the case of degenerate meromorphic mappings, in \cite{RW}, Ru and Wang gave a second main theorem for moving hyperplanes with counting function truncated to level $n$. And then, the result of Ru-Wang was improved by Thai-Quang \cite{TQ08} and Quang-An \cite{QA}. In 2016, the author have improved and extended all those results to a better second main theorem. To state their results, we recall the following. 

Let $a_1,\dots,a_q$ $(q \geq n+1)$ be $q$ meromorphic mappings of $\C^m$ into the dual space $\P^n(\C)^*$ with reduced representations $a_i = (a_{i0}: \dots : a_{in})\ (1\le i \le q).$ We say that $a_1,\dots,a_q$ are located in general position if $\det (a_{i_kl}) \not \equiv 0$ for any $1\le i_0<i_1<\cdots <i_n\le q.$ Let $\mathcal {M}_m$ be the field of all meromorphic functions on $\C^m$. Denote by $\mathcal R_{\{a_i\}_{i=1}^q} \subset \mathcal {M}_m$ the smallest subfield which contains $\C$ and all $ \frac {a_{ik}}{a_{il}}\text { with } a_{il}\not\equiv 0.$ Thoughout this paper, if without any notification, the notation $\mathcal R$ is always stands for $\mathcal R_{\{a_i\}_{i=1}^q}$.

In 2004, M. Ru and J. Wang proved the following.

\vskip0.2cm
\noindent{\bf Theorem A }\cite[Theorem 1.3]{RW} {\it Let $f: \mathbb C\to\mathbb P^n(\mathbb C)$ be a holomorphic map. Let $\{a_j\}_{j=1}^q$ be moving hyperplanes of $\P^n(\C)$ in general position such that $(f,a_j)\not\equiv0$\ $(1\leq j\leq q)$. If $q\ge 2n+1$ then
$$\big|\big|\ \displaystyle\frac{q}{n(2n+1)}T_f(r)\leq\sum_{i=1}^q N^{[n]}_{(f_i,a)}(r)+o(T_f(r))+O(\max\limits_{1\leq i\leq q}T_{a_i}(r)).$$}

\vskip0.2cm
Here, by the notation $``|| \ P"$  we mean the assertion $P$ holds for all $r \in [0,\infty)$ outside a Borel subset $E$ of the interval $[0,\infty)$ with $\int_E dr<\infty$. 

In 2008, D. D. Thai and S. D. Quang improved the above result by increasing the coefficent $\frac{q}{n(2n+1)}$ in front of the characteristic function to $\frac{q}{2n+1}$. In 2016, S D. Quang \cite{Q14} improved these result to the following.

\vskip0.2cm
\noindent{\bf Theorem B }\cite[Theorem 1.1]{Q14} {\it Let $f: \mathbb C^m\to\mathbb P^n(\mathbb C)$ be a meromorphic mapping. Let $\{a_j\}_{j=1}^q$\ $(q\geq 2n-k+2)$ be meromorphic mappings of $\mathbb C^m$ into $\mathbb P^n(\mathbb C)^*$ in general position such that $(f,a_j)\not\equiv0$\ $(1\leq j\leq q)$, where $\rank_{R\{a_j\}}(f)=k+1$. Then the following assertion holds:\\
(a) $\big|\big|\ \displaystyle\frac{q}{2n-k+2}T_f(r)\leq\sum_{i=1}^q N^{[k]}_{(f_i,a)}(r)+o(T_f(r))+O(\max\limits_{1\leq i\leq q}T_{a_i}(r)),$\\
(b) $\big|\big|\ \displaystyle\frac{q-(n+2k-1)}{n+k+1}T_f(r)\leq\sum_{i=1}^q N^{[k]}_{(f_i,a)}(r)+o(T_f(r))+O(\max\limits_{1\leq i\leq q}T_{a_i}(r)).$}
\vskip0.2cm

The main purpose of the present paper is to establish a stronger second main theorem for  meromorphic mappings of $\C^m$ into $\P^n(\C)$ and moving hyperplanes. Namely, we will prove the following theorem.
\begin{theorem}\label{1.1}
Let $f :\C^m \to \P^n(\C)$ be a meromorphic mapping. Let $\{a_i\}_{i=1}^q \ (q\ge 2n-k+2)$ be meromorphic mappings of $\C^m$ into $\P^n(\C)^*$ in general position such that $(f,a_i)\not\equiv 0\ (1\le i\le q),$ where $k+1=\rank_{\mathcal R\{a_i\}}(f)$. Then the following assertions hold:
\begin{align*}
&\text{(a)}\ || \ \dfrac {q-(n-k)}{n+2}T_f(r) \le \sum_{i=1}^q N_{(f,a_i)}^{[k]}(r) + o(T_f(r)) + O(\max_{1\le i \le q}T_{a_i}(r)),\\
&\text{(b)}\ || \ \dfrac{q-2(n-k)}{k(k+2)}T_f(r) \le \sum_{i=1}^q N_{(f,a_i)}^{[1]}(r) + o(T_f(r)) + O(\max_{1\le i \le q}T_{a_i}(r)).
\end{align*}
\end{theorem}
Here by $\rank_{\mathcal R\{a_i\}}(f)$ we denote the rank of the set $\{f_0,f_1,\ldots, f_n\}$ over the field $\mathcal R\{a_i\}$, where $(f_0:f_1:\cdots :f_n)$ is a representation of the mapping $f$.

\textit{\underline{Remark:}} 1) The assertion (a) is an improvement of Theorem B.

2) It is easy to see that $\frac{q-2(n-k)}{k(k+2)}\ge \frac{q}{n(n+2)}$. Therefore, the assertion (b) immediately implies the following corollary.
\begin{corollary}\label{1.3}
With the assumptions of Theorem A, we have
$$|| \dfrac{q}{n(n+2)}T_f(r) \le \sum_{i=1}^q N_{(f,a_i)}^{[1]}(r) + o(T_f(r)) + O(\max_{1\le i \le q}T_{a_i}(r)).$$
\end{corollary}

In order to prove the above result, beside developing the method used in \cite{Q14,RW,TQ05}, we also propose some new techniques. Firstly, we will rearrange the family hyperplanes in the increasing order of the values of the counting functions (of their inverse images). After that, we find the smallest number of the first hyperplanes in this order such that the sum of their counting functions exceed the characteristic functions. And then, we have to compare the characteristic functions with this sum of counting functions with explicitly estimating the truncation level. From that, we deduce the second main theorem. 

For the case where the number of moving hyperplanes is large enough, we will prove a better second main theorem as follows.
\begin{theorem}\label{1.2}
With the assumptions of Theorem \ref{1.1}, we assume further more that $q\ge (n-k)(k+1)+n+2$. Then we have
\begin{align*}
|| \ \dfrac {q}{k+2}T_f(r) \le \sum_{i=1}^q N_{(f,a_i)}^{[k]}(r) + o(T_f(r)) + O(\max_{1\le i \le q}T_{a_i}(r)).
\end{align*}
\end{theorem}

In this case, we may see that the coefficient in front of the characteristic functions are exactly the same as the case where the mappings are assumed to be non-degenerate.

{\bf Acknowledgements.} This work was done during a stay of the author at the Vietnam Institute for Advanced Study in Mathematics. He would like to thank the institute for their support. This research is funded by Vietnam National Foundation for Science and Technology Development (NAFOSTED) under grant number 101.04-2018.01.

\section{Basic notions and auxiliary results from Nevanlinna theory}

Throughout this paper, we use the standart notation on Nevanlina theory from \cite{Q14,Q15} and \cite{TQ08}. For a meromorphic mapping $f:\C^m\to\P^n(\C)$, we denote by $T_f(r)$ its characteristic funtion. 
Let $\varphi$ be a meromorphic funtion on $\C^m$. We denote by $\nu_\varphi$ its divisor, $N^{[k]}_\varphi (r)$ the counting function with the trucation level $k$ of its zeros divisor and $m(r,\varphi)$ its proximity function. The lemma on logarithmic derivative in Nevanlinna theory is stated as follows.
\begin{lemma}[{\cite[Lemma 3.11]{Shi}}]\label{2.1}
Let $f$ be a nonzero meromorphic function on $\C^m.$ Then 
$$\biggl|\biggl|\quad m\biggl(r,\dfrac{\mathcal{D}^\alpha (f)}{f}\biggl)=O(\log^+T_f(r))\ (\alpha\in \Z^m_+).$$
\end{lemma}
The first main theorem states that
$$ T_{\varphi}(r)=m(r,\varphi) +N_{\frac{1}{\varphi}}(r).$$

We assume that thoughout this paper, the homogeneous coordinates of $\P^n(\C)$ is chosen so that for each given meromorphic mapping $a=(a_0:\cdots :a_n)$ of $\C^m$ into $\P^n(\C)^*$ then $a_{0}\not\equiv 0$. We set
$$ \tilde a_i=\dfrac{a_i}{a_0}\text{ and }\tilde a=(\tilde a_0:\tilde a_1:\cdots:\tilde a_n).$$
Supposing that $f$ has a reduced representation $f=(f_0:\cdots :f_n),$ we put $(f,a):=\sum_{i=0}^{n}f_ia_{i}$ and $(f,\tilde a):=\sum_{i=0}^{n}f_i\tilde a_{i}.$

Let $\{a_i\}_{i=1}^q$ be $q$ meromorphic mappings of $\C^m$ into $\P^n(\C)^*$  with reduced representations $a_i=(a_{i0}:\cdots :a_{in})\ (1\le i\le q).$ 
\begin{definition}
The family $\{a_i\}_{i=1}^q$ is said to be in general position  if  
$$\det (a_{i_jl};0\le j\le n,0\le l\le n) \not\equiv 0$$
for any $1\le i_0\le\cdots\le i_n\le q$.
\end{definition}

\begin{definition}
A subset $\mathcal {L}$ of $\mathcal {M}$ (or $\mathcal {M}^{n+1}$) is said to be minimal over the field $\mathcal R$ if it is linearly dependent over $\mathcal {R}$ and each proper subset of $\mathcal L$ is linearly independent over  $\mathcal {R}.$
\end{definition}

Repeating the argument in \cite[Proposition 4.5]{Fu}, we have the following proposition.
\begin{proposition}[{see \cite[Proposition 4.5]{Fu}}]\label{2.2}
Let $\Phi_0,\ldots,\Phi_k$ be meromorphic functions on $\C^m$ such that $\{\Phi_0,\ldots,\Phi_k\}$ 
are  linearly independent over $\C.$
Then  there exists an admissible set $\{\alpha_i=(\alpha_{i1},\ldots,\alpha_{im})\}_{i=0}^k \subset \Z^m_+$ with $|\alpha_i|=\sum_{j=1}^{m}|\alpha_{ij}|\le k \ (0\le i \le k)$ satisfying the following two properties:

(i)\  $\{{\mathcal D}^{\alpha_i}\Phi_0,\ldots,{\mathcal D}^{\alpha_i}\Phi_k\}_{i=0}^{k}$ is linearly independent over $\mathcal M,$\ i.e., \ 
$\det{({\mathcal D}^{\alpha_i}\Phi_j)}\not\equiv 0,$ 

(ii) $\det \bigl({\mathcal D}^{\alpha_i}(h\Phi_j)\bigl)=h^{k+1}\det \bigl({\mathcal D}^{\alpha_i}\Phi_j\bigl)$ for
any nonzero meromorphic function $h$ on $\C^m.$
\end{proposition}

\section{Proof of Theorem \ref{1.1} and Theorem \ref{1.2}}

In order to prove Theorem \ref{1.1} we need the following lemma, which is an improvement of Lemma 3.1 in \cite{Q14}.
\begin{lem}\label{3.1}
Let  $f:\C^m\rightarrow\P^n(\C)$ be a meromorphic mapping. Let $\{a_i\}_{i=1}^p$ be $p$ meromorphic mappings of $\C^m$ into $\P^n(\C)^*$ in general position with $\rank\{(f,\tilde a_i);1\le i\le q\}=\rank_{\mathcal R}(f)$, where $\mathcal R=\mathcal R\{a_i\}_{i=1}^p$. Assume that there exists a partition $\{1,\ldots,q\}=I_1\cup I_2\cdots\cup I_l$ satisfying:\\
$\mathrm{(i)}$ \ $\{(f,\tilde a_i)\}_{i\in I_1}$ is minimal over $\mathcal R$, $\{(f,\tilde a_i)\}_{i\in I_t}$ is linearly independent over $\mathcal {R}\ (2\le t \le l), $ \\
$\mathrm{(ii)}$ \ For any $2\le t\le l,i\in I_t,$ there exist meromorphic functions $c_i\in \mathcal {R}\setminus\{0\}$ such that 
$$\sum_{i\in I_t}c_i(f,\tilde a_i)\in \biggl(\bigcup_{j=1}^{t-1}\bigcup_{i\in I_j}(f,\tilde a_i) \biggl)_{\mathcal {R}}.$$
Then we have
$$ T_f(r)\le \sum_{i=1}^{t}\sum_{j\in I_i}^qN^{[n_i]}_{(f,a_j)}(r)+ o(T_f(r)) + O(\max_{1\le i \le p}T_{a_i}(r)),$$
where $n_1=\sharp I_1-2$ and $n_t=\sharp I_t-1$ for $t=2,...,l$.
\end{lem}

\textbf{Proof.}\ Let $f=(f_0:\cdots :f_n)$ be a reduced representation of $f$. By changing the homogeneous coordinate system of $\P^n(\C)$ if necessary, we may assume that $f_0\not\equiv 0.$ 
Without loss of generality, we may assume that $I_1=\{1,\ldots.,k_1\}$ and
$$I_t=\{k_{t-1}+1,\ldots, k_t\}\ (2\le t \le l),\text{ where }1=k_0<\cdots< k_l=q.$$ 

Since $\{(f,\tilde a_i)\}_{i\in I_1}$ is minimal over $\mathcal R$, there exist $c_{1i}\in\mathcal {R}\setminus \{0\}$ such that 
$$\sum_{i=1}^{k_1}c_{1i}\cdot (f, \tilde a_i)=0.$$
Define $c_{1i}=0$ for all $i>k_1.$ Then 
$$\sum_{i=1}^{k_l} c_{1i}\cdot (f,\tilde a_i)=0.$$
Because ${\{c_{1i}(f,\tilde a_i)\}}_{i=k_0+1}^{k_1}$ is linearly independent over $\mathcal R,$ Proposition \ref{2.2} yields that there exists an admissible set $\{\alpha_{1(k_0+1)},\ldots,\alpha_{1k_1}\}\subset \Z^m_+$ \ $(|\alpha_{1i}|\le k_1-k_0-1=n_1)$ such that the matrix
$$\ A_1=\left (\mathcal {D}^{\alpha_{1i}}(c_{1j}(f,\tilde a_j));k_0+1\le i,j\le k_1 \right)$$
has nonzero determinant.

Now consider $t\ge 2.$ 
By the construction of the set $I_t$, there exist meromorphic mappings $c_{ti}\not\equiv 0\ (k_{t-1}+1\le i\le k_t)$ such that 
$$\sum_{i=k_{t-1}+1}^{k_t}c_{ti}\cdot (f,\tilde a_i)\in 
\biggl(\bigcup_{j=1}^{t-1}\bigcup_{i\in I_t}{(f,\tilde a_i)}\biggl)_{\mathcal {R}}.$$
Therefore, there exist meromorphic mappings $c_{ti}\in \mathcal {R}\ (1\le i\le k_{t-1})$ 
such that 
$$\sum_{i=1}^{k_t}c_{ti}\cdot (f,\tilde a_i)=0.$$
Define $c_{ti}=0$ for all $i>k_t.$ Then 
$$\sum_{i=1}^{k_l}c_{ti}\cdot (f,\tilde a_i)=0.$$
Since $\{c_{ti}(f,\tilde a_i)\}_{i=k_{t-1}+1}^{k_t}$ is $\mathcal {R}$-linearly independent, by again Proposition \ref{2.2} there exists an admissible set $\{\alpha_{t(k_{t-1}+1)},\ldots,\alpha_{tk_t}\}\subset \Z^m_+$ \ $(|\alpha_{ti}|\le k_t-k_{t-1}-1=n_t)$ such that the matrix
$$\ A_t=\left (\mathcal {D}^{\alpha_{ti}}(c_{1j}(f,\tilde a_j));k_{t-1}+1\le i,j\le k_t \right)$$
has nonzero determinant.

Consider the following $(k_l-1)\times k_l$ matrix 
\begin{align*}
T&=\left [\biggl (\mathcal {D}^{\alpha_{ti}}(c_{tj}(f,\tilde a_j));1\le t\le l, k_{t-1}+1\le i\le k_t\biggl );1\le j\le k_l \right]
\end{align*}
Denote by $D_i$ the subsquare matrix obtained by deleting the $i$-th column 
of the minor  matrix  $T$. Since the sum of each row of $T$ is zero, we have
$$\det D_i={(-1)}^{i-1}\det D_1={(-1)}^{i-1}\prod_{j=1}^{l}\det A_j.$$

Since $\{a_i\}_{i=1}^q$ is in general position, we have 
$$\det (\tilde a_{ij}, \ 1\le i\le n+1,0\le j\le n )\not\equiv 0.$$
By solving the linear equation system $(f,\tilde a_i)=\tilde a_{i0}\cdot f_0+\ldots +\tilde a_{in}\cdot f_n \ (1\le i\le n+1),$
we obtain 
\begin{align}\label{+}f_v=\sum_{i=1}^{n+1}A_{vi}(f,\tilde a_{i})\ (A_{vi}\in\mathcal R)\text{ for each }0\le v \le n.
\end{align}
Put 
$\Psi(z)=\sum_{i=1}^{n+1}\sum_{v=0}^n |A_{vi}(z)|\ (z\in \C^m).$
Then
$$\ \ ||f(z)||\le \Psi (z)\cdot \max_{1\le i\le n+1}\bigl (|(f,\tilde a_i)(z)|\bigl )\le \Psi (z)\cdot \max_{1\le i\le q}\bigl (|(f,\tilde a_i)(z)|\bigl )\ (z\in \C^m),$$
and
\begin{align*}
\int\limits_{S(r)} \log^+\Psi (z) \eta \le \sum_{i=1}^{n+1}\sum_{v=0}^n T(r,A_{vi}) +O(1)=  O(\max_{1\le i\le q}T_{a_i}(r))+O(1).
\end{align*}

Fix $z_0 \in \C^m\setminus\bigcup_{j=1}^q\biggl (\supp (\nu^0_{(f,\tilde a_j)})\cup \supp (\nu^\infty_{(f,\tilde a_j)})\biggl ).$ Take $i\ (1\le i \le q)$ such that
$$|(f,\tilde a_i)(z_0)|=\max_{1\le j\le q}|(f,\tilde a_j)(z_0)|.$$
 Then
\begin{align*}
\log\dfrac{|\det D_1(z_0)|.||f(z_0)||}{\prod_{j=1}^{q}|(f,\tilde a_j)(z_0)|}& \le\log^+\biggl ( \Psi (z_0)\cdot \biggl(\dfrac{|\det D_i(z_0)|}{\prod_{j=1,j\ne i}^{q}|(f,\tilde a_j)(z_0)|}\biggl)\biggl  )\\
& \le\log^+\biggl(\dfrac{|\det D_i(z_0)|}{\prod_{j=1,j\ne i}^{q}|(f,\tilde a_j)(z_0)|}\biggl)+\log^+\Psi (z_0).
\end{align*}

Thus, for each $z\in \C^m\setminus\bigcup_{j=1}^q\biggl (\supp (\nu^0_{(f,\tilde a_j)})\cup \supp (\nu^\infty_{(f,\tilde a_j)})\biggl ),$ we have
\begin{align*}
\log\dfrac{|\det D_1(z)|.||f(z)||}{\prod_{i=1}^{q}|(f,\tilde a_i)(z)|} \le\sum_{i=1}^{q}\log^+\biggl(\dfrac{|\det D_i(z)|}{\prod_{j=1,j\ne i}^{q}|(f,\tilde a_j)(z)|}\biggl)+\log^+ \Psi (z).
\end{align*}
Hence
\begin{align}
\log ||f(z)||\le \log \dfrac{\prod_{i=1}^{q}|(f,\tilde a_i)(z)|}{|\det D_1(z)|}+\sum_{i=1}^{q}\log^+\biggl(\dfrac{|\det D_i(z)|}{\prod_{j=1,j\ne i}^{q}|(f,\tilde a_j)(z)|}\biggl)+\log^+ \Psi (z).
\end{align}
 
Integrating both sides of the above inequality and using Jensen's formula and the lemma on logarithmic derivative, we have
\begin{align}\nonumber
||\ \ T_f(r)\le &N_{\prod_{i=1}^{q}(f, \tilde a_i)}(r)-N(r,\nu_{\det D_1})+ \sum_{i=1}^{q} m\biggl(r,\dfrac{\det D_i}{\prod_{j=1,j\ne i}^{q}(f,\tilde a_j)}\biggl)
+O(\max_{1\le i \le q}T_{a_i}(r))\\
\label{3.3}
=&N_{\prod_{i=1}^{q}(f, \tilde a_i)}(r)-N(r,\nu_{\det D_1})+ O(\log^+T_f(r))+O(\max_{0\le i \le q-1}T_{a_i}(r)).
\end{align}

\begin{claim}
$||\ N_{\prod_{i=1}^{q}(f, \tilde a_i)}(r)-N(r,\nu_{\det D_1})\le \sum_{s=1}^l\sum_{i\in I_s}^qN^{[n_s]}_{(f,a_i)}(r)+O(\max_{1\le i\le q}T_{a_i}(r)).$
\end{claim}
Indeed, fix $z\in\C^m\setminus I(f)$, where $I(f)=\{f_0=\cdots f_n=0\}$. We call $i_0$ the index satisfying 
$$\nu^0_{(f,\tilde a_{i_0})}(z)=\min_{1\le i\le n+1}\nu^0_{(f,\tilde a_i)}(z).$$

For each $i\ne i_0, i\in I_s$, we easily have
$$ \nu_{\mathcal {D}^{\alpha_{sk_{s-1}+j}}(c_{si}(f,\tilde a_{i}))}(z)\ge \max\{0,\nu_{(f\tilde a_i)}^{0}(z)-n_s\}-C\bigl (2\nu_{c_{si}}^{\infty}(z)+\nu^0_{a_{i0}}(z)\bigl ),$$
where $C$ is a fixed constant. 

Since each element of the matrix $D_{i_0}$ is of the form $\mathcal {D}^{\alpha_{sk_{s-1}+j}}(c_{si}(f,\tilde a_{i}))\ (i\ne i_0)$, one estimates
\begin{align}\label{3.4}
\nu_{D_1}(z)=\nu_{D_{i_0}}(z)
\ge\sum_{s=1}^l\sum_{\overset{i\in I_s}{i\ne i_0}}\left (\max\{0,\nu_{(f\tilde a_i)}^{0}(z)-n_s\}-(k+1)\bigl (2\nu_{c_{si}}^{\infty}(z)+\nu^0_{a_{i0}}(z)\bigl )\right ).
\end{align}
We see that there exists $v_0\in\{0,\ldots,n\}$ with $f_{v_0}(z)\ne 0$. Then by (\ref{+}), there exists $i_1\in\{1,\ldots,n+1\}$ such that $A_{v_0i_1}(z)\cdot (f,\tilde a_{i_1})(z)\ne 0$. Thus
\begin{align}\label{3.5}
\nu^0_{(f,\tilde a_{i_0})}(z)\le \nu^0_{(f,\tilde a_{i_1})}(z)\le\nu^\infty_{A_{v_0i_1}}(z)\le\sum_{A_{vi}\not\equiv 0}\nu^\infty_{A_{vi}}(z).
\end{align}
Combining the inequalities (\ref{3.4}) and (\ref{3.5}), we have
\begin{align*}
\nu^{0}_{\prod_{i=1}^{q}(f,\tilde a_i)}&(z)-\nu_{\det D_1}(z)\\
&\le\sum_{s=1}^l\sum_{\overset{i\in I_s}{i\ne i_0}}\left (\min\{\nu_{(f,\tilde a_i)}^{0}(z),n_s\}+(k+1)\bigl (2\nu_{c_{si}}^{\infty}(z)+\nu^0_{a_{i0}}(z)\bigl )\right )+\sum_{A_{vi}\not\equiv 0}\nu^\infty_{A_{vi}}(z)\\
&\le\sum_{s=1}^l\sum_{i\in I_s}\left (\min\{\nu_{(f,\tilde a_i)}^{0}(z),n_s\}+(k+1)\bigl (2\nu_{c_{si}}^{\infty}(z)+\nu^0_{a_{i0}}(z)\bigl )\right )+\sum_{A_{vi}\not\equiv 0}\nu^\infty_{A_{vi}}(z).
\end{align*}
Integrating both sides of this inequality, we easily obtain
\begin{align}\label{3.6}
||\ \ N_{\prod_{i=1}^{q}(f, \tilde a_i)}(r)-N(r,\nu_{\det D_1})\le\sum_{s=1}^l\sum_{i\in I_s}N^{[n_s]}_{(f,a_i)}(r)+O(\max_{1\le i\le q}T_{a_i}(r)).
\end{align}
The claim is proved.

From the inequalities (\ref{3.3}) and the claim, we get
$$ ||\ \ T_f(r)\le \sum_{s=1}^l\sum_{i\in I_s}N^{[n_s]}_{(f,a_i)}(r)+O(\log^+T_f(r))+O(\max_{1\le i\le q}T_{a_i}(r)). $$
The lemma is proved.\hfill$\square$

\vskip0.2cm
\noindent
{\bf Proof of Theorem \ref{1.1}.}

We denote by $\mathcal I$ the set of all permutations of $q$-tuple $(1,\ldots,q)$. For each element $I=(i_1,\ldots,i_q)\in\mathcal I$, we set
\begin{align*}
N_I&=\{r\in\R^+;N^{[k]}_{(f,a_{i_1})}(r)\le\cdots\le N^{[k]}_{(f,a_{i_q})}(r)\},\\
M_I&=\{r\in\R^+;N^{[1]}_{(f,a_{i_1})}(r)\le\cdots\le N^{[1]}_{(f,a_{i_q})}(r)\}.
\end{align*}
We now consider an element $I=(i_1,\ldots,i_q)$ of $\mathcal I$. We will construct subsets $I_t$ of the set $A_1=\{1,\ldots,{2n-k+2}\}$ as follows.

We choose a subset $I_1$ of $A$ which is the minimal subset of $A$ satisfying that $\{(f,\tilde a_{i_j})\}_{j\in I_1}$ is minimal over $\mathcal R$. If $\rank_{\mathcal R}\{(f,\tilde a_{i_j})\}_{j\in I_1}= k+1$ then we stop the process. 

Otherwise, set $I_1'=\{i; (f,\tilde a_i)\in \bigl (\{(f,\tilde a_{i_j})\}_{j\in I_1}\bigl )\}$, $A_2=A_1\setminus (I_1\cup I_1')$ and see that $\sharp I_1\cup I_1'\le n+1$. We consider the following two cases:
\begin{itemize}
\item Case 1. Suppose that $\sharp A_2\ge n+1$. Since $\{\tilde a_{i_j}\}_{j\in A_2}$ is in general position, we have
$$ \left ((f,\tilde a_{i_j}); j\in A_2\right )_{\mathcal R}=\left (f_0,\ldots,f_n\right )_{\mathcal R}\supset \left ((f,\tilde a_{i_j}); j\in I_1\right )_{\mathcal R}\not\equiv 0.$$
\item  Case 2. Suppose that $\sharp A_2< n+1$. Then we have the following:
\begin{align*}
&\dim_{\mathcal R}\left ((f,\tilde a_{i_j}); j\in I_1\right )_{\mathcal R}\ge k+1-(n+1-\sharp I_1\cup I_1')=k-n+\sharp I_1\cup I_1',\\ 
& \dim_{\mathcal R}\left ((f,\tilde a_{i_j}); j\in A_2\right )_{\mathcal R}\ge k+1-(n+1-\sharp A_2)=k-n+\sharp A_2.
\end{align*}
We note that $\sharp I_1\cup I_1'+\sharp A_2=2n-k+2$. Hence the above inequalities imply that
\begin{align*}
\dim_{\mathcal R}&\biggl (\bigl ((f,\tilde a_{i_j}); j\in I_1\bigl )_{\mathcal R}\cap\bigl ((f,\tilde a_{i_j}); j\in A_2\bigl )_{\mathcal R}\biggl )\\
&\ge\dim_{\mathcal R}\left ((f,\tilde a_{i_j}); j\in I_1\cup I_1'\right )_{\mathcal R}+\dim_{\mathcal R}\left ((f,\tilde a_{i_j}); j\in A_2\right )_{\mathcal R}-(k+1)\\
&=k-n+\sharp I_1\cup I_1'+k-n+\sharp A_2-(k+1)=1.
\end{align*}
\end{itemize}
Therefore, from the above two cases, we see that
$$ \bigl ((f,\tilde a_{i_j}); j\in I_1\bigl )_{\mathcal R}\cap\bigl ((f,\tilde a_{i_j}); j\in A_2\bigl )_{\mathcal R}\ne \{0\}. $$
Therefore, we may chose a subset $I_2\subset A_2$ which is the minimal subset of $A_2$ satisfying that there exist nonzero meromorphic functions $c_i\in\mathcal R\ (i\in I_2)$,
$$\sum_{i\in I_2}c_i(f,\tilde a_i)\in \biggl(\bigcup_{i\in I_1}(f,\tilde a_i) \biggl)_{\mathcal {R}}.$$
We see that $\sharp I_2\ge 2$. By the minimality of the set $I_2$, the family $\{(f,\tilde a_{i_j})\}_{j\in I_2}$ is linearly independent over $\mathcal R$, and hence $\sharp I_2\le k+1$ and 
$$\sharp (I_2\cup I_2)\le\min\{2n-k+2, n+k+1\}.$$ 

Moreover, we will show that 
$$\dim\biggl (\bigl ((f,\tilde a_{i_j}); j\in I_1\bigl )_{\mathcal R}\cap\bigl ((f,\tilde a_{i_j}); j\in A_2\bigl )_{\mathcal R}\biggl )=1.$$
Indeed, suppose contrarily there exist two linearly independent vectors $x,y\in \bigl ((f,\tilde a_{i_j}); j\in I_1\bigl )_{\mathcal R}\cap\bigl ((f,\tilde a_{i_j}); j\in I_2\bigl )_{\mathcal R}$, with 
\begin{align*}
x&=\sum_{i\in I_2}x_i(f,\tilde a_i)\in \bigl ((f,\tilde a_{i_j}); j\in I_1\bigl )_{\mathcal R},\\ 
y&=\sum_{i\in I_2}y_i(f,\tilde a_i)\in \bigl ((f,\tilde a_{i_j}); j\in I_1\bigl )_{\mathcal R},
\end{align*}
where $x_i,y_i\in\mathcal R$. By the minimality of the se $I_2$, all functions $x_i,y_i$ are not zero. Therefore, fixing $i_0\in I_2$, we have
$$ \sum_{\overset{i\in I_2}{i\ne i_0}}(y_0x_i-x_0y_i)(f,\tilde a_i)\in \bigl ((f,\tilde a_{i_j}); j\in I_1\bigl )_{\mathcal R}. $$
Since $x,y$ are linearly independent, the left hand side is not zero. This contradics the minimality of the set $I_2$. Hence 
$$\dim\biggl (\bigl ((f,\tilde a_{i_j}); j\in I_1\bigl )_{\mathcal R}\cap\bigl ((f,\tilde a_{i_j}); j\in I_2\bigl )_{\mathcal R}\biggl )=1.$$

On the other hand, we will see that $\sharp I_1\cup I_2\le n+2$. If 
$\rank_{\mathcal R}\{(f,\tilde a_{i_j})\}_{j\in I_1\cup I_2}= k+1$ then we stop the process.

Otherwise, by repeating the above argument, we have a subset $I_2'=\{i; (f,\tilde a_i)\in \bigl (\{(f,\tilde a_{i_j})\}_{j\in I_1\cup I_2}\bigl )\}$, a subset $I_3$ of $A_3=A_1\setminus (I_1\cup I_2\cup I_2')$, which satisfy the following:
\begin{itemize}
\item there exist nonzero meromorphic functions $c_i\in\mathcal R\ (i\in I_3)$ so that
$$\sum_{i\in I_3}c_i(f,\tilde a_i)\in \biggl(\bigcup_{i\in I_1\cup I_2}(f,\tilde a_i) \biggl)_{\mathcal {R}},$$
\item $\{(f,\tilde a_{i_j})\}_{j\in I_3}$ is linearly independent over $\mathcal R$,
\item $2\le \sharp I_3\le k+1$ and  $\sharp (I_1\cup\cdots\cup I_3)\le \min\{2n-k+2, n+k+1\}$,
\item $\dim\biggl (\bigl ((f,\tilde a_{i_j}); j\in I_1\cup I_2\bigl )_{\mathcal R}\cap\bigl ((f,\tilde a_{i_j}); j\in I_3\bigl )_{\mathcal R}\biggl )=1.$
\end{itemize}

Continuing this process, we get a sequence of subsets $I_1,\ldots,I_l$, which satisfy: 
\begin{itemize}
\item[(1)] $\{(f,\tilde a_{i_j})\}_{j\in I_1}$ is minimal over $\mathcal R$, $\sharp I_t\ge 2$ and $\{(f,\tilde a_{i_j})\}_{j\in I_t}$ is linearly independent over $\mathcal {R}\ (2\le t \le l), $
\item[(2)] for any $2\le t\le l, j\in I_t,$ there exist meromorphic functions $c_j\in \mathcal {R}\setminus\{0\}$ such that 
$$\sum_{j\in I_t}c_j(f,\tilde a_{i_j})\in \biggl(\bigcup_{s=1}^{t-1}\bigcup_{j\in I_s}(f,\tilde a_{i_j}) \biggl)_{\mathcal {R}},$$\\
$$\text{ and }\dim\biggl (\bigl ((f,\tilde a_{i_j}); j\in I_1\cup\cdots\cup I_{t-1}\bigl )_{\mathcal R}\cap\bigl ((f,\tilde a_{i_j}); j\in I_t\bigl )_{\mathcal R}\biggl )=1,$$
\item[(3)] $\rank_{\mathcal R}\{(f,\tilde a_{i_j})\}_{j\in  I_1\cup\cdots\cup I_{l}}= k+1$.
\end{itemize}
If $\sharp I_1=2$ we will remove one element from $I_1$ and combine the remaining element with $I_2$ to become a new set $I_1$. Therefore, we will get a sequence $I_1,...,I_l$ which satisfy the above three properties and $\sharp I_1\ge 3, \sharp I_t\ge 2\ (2\le t\le l)$. We set $n_1=\sharp I_1-2, n_s=\sharp I_s-1\ (2\le s\le l), n_0=\max_{1\le s\le l}n_s, J=I_1\cup\cdots\cup I_l$ and $d+2=\sharp J$. One estimates
\begin{align*}
&(n_1+2)+(n_2+1)+\cdots +(n_l+1)=d+2,\\ 
&(n_1+1)+n_2+\cdots + n_l=k+1.
\end{align*}
Since the rank of the set of any $n+1$ functions $(f,\tilde a_i)'$s is equal to $k+1$, we have
$$ (n+1)-\sharp (I_1\cup\cdots \cup I_{l-1})\ge (k+1)-\rank\{(f,\tilde a_i);\ i\in I_1\cup\cdots \cup I_{l-1}\}, $$
$$ \mathrm{i.e., }(n+1)- (n_1+2)-(n_2+1)-\cdots -(n_{l-1}+1)\ge (k+1)-(n_1+1)-n_2-\cdots -n_{l-1}.$$
This implies that
$$ d+2\le n+2. $$
On the other hand, we see that $k+1+l=d+2$, and hence
\begin{align*}
n_s=k-\sum_{\overset{i=1}{i\ne s}}^ln_i\le k-(l-1)\le\frac{k(k+2)}{k+l+1}=\frac{k(k+2)}{d+2}. 
\end{align*}
Thus $n_0\le\frac{k(k+2)}{d+2}$.

Now the family of subsets $I_1,\ldots,I_l$ satisfies the assumptions of the Lemma \ref{3.1}. Therefore, we have
\begin{align}\label{new1}
||\ T_f(r)\le\sum_{s=1}^l\sum_{j\in I_s}N^{[n_s]}_{(f,a_{j})}+o(T_f(r))+O(\max_{1\le i \le q}T_{a_i}(r)).
\end{align}

(a) For all $r\in N_I$ (may be outside a finite Borel measure subset of $\R^+$), from (\ref{new1}) we have
\begin{align*}
||\ T_f(r)&\le\sum_{j\in J}N^{[k]}_{(f,a_{j})}+o(T_f(r))+O(\max_{1\le i \le q}T_{a_i}(r))\\
&\le\dfrac{\sharp J}{q-(2n-k+2)+\sharp J}\biggl (\sum_{j\in J}N^{[k]}_{(f,a_{i_j})}(r)+\sum_{j=2n-k+3}^qN^{[k]}_{(f,a_{i_j})}(r)\biggl )\\
&+o(T_f(r)) + O(\max_{1\le i \le q}T_{a_i}(r)).
\end{align*}
Since $\sharp J=d+2\le n+2$, the above inequality implies that 
\begin{align}\label{3.8}
||\ T_f(r)\le\dfrac{n+2}{q-(n-k)}\sum_{i=1}^qN^{[k]}_{(f,a_{i})}(r)+o(T_f(r))+O(\max_{1\le i \le q}T_{a_i}(r)),\quad r\in N_I.
\end{align}

We see that $\bigcup_{I\in\mathcal I}N_I=\R^+$ and the inequality (\ref{3.8}) holds for every $r\in N_I, I\in\mathcal I$. This yields that
$$ T_f(r)\le\dfrac{n+2}{q-(n-k)}\sum_{i=1}^qN^{[k]}_{(f,a_{i})}(r)+o(T_f(r))+O(\max_{1\le i \le q}T_{a_i}(r)) $$ 
for all $r$ outside a finite Borel measure subset of $\R^+$. Thus
$$ ||\ \dfrac{q-(n-k)}{n+2}T_f(r)\le\sum_{i=1}^qN^{[k]}_{(f,a_{i})}(r)+o(T_f(r))+O(\max_{1\le i \le q}T_{a_i}(r)). $$
The assertion (a) is proved. 

(b) We repeat the same argument as in the proof of the assertion (a). For all $r\in M_I$ (may be outside a finite Borel measure subset of $\R^+$) we have
\begin{align*}
||\ T_f(r)&\le\sum_{s=1}^l\sum_{j\in I_s}N^{[n_s]}_{(f,a_{j})}+o(T_f(r))+O(\max_{1\le i \le q}T_{a_i}(r))\\
&\le\sum_{j\in J}n_0N^{[1]}_{(f,a_{j})}+o(T_f(r))+O(\max_{1\le i \le q}T_{a_i}(r))\\
&\le n_0\cdot\frac{d+2}{q-(2n-k+2)+d+2}\biggl (\sum_{j\in J}N^{[1]}_{(f,a_{i_j})}(r)+\sum_{j=2n-k+3}^qN^{[1]}_{(f,a_{i_j})}(r)\biggl )\\
&+o(T_f(r)) + O(\max_{1\le i \le q}T_{a_i}(r))\\
&\le \frac{k(k+2)}{q-(2n-k+2)+d+2}\sum_{i=1}^qN^{[1]}_{(f,a_{i})}(r)+o(T_f(r))+O(\max_{1\le i \le q}T_{a_i}(r))\\
&\le \frac{k(k+2)}{q-2(n-k)}\sum_{i=1}^qN^{[1]}_{(f,a_{i})}(r)+o(T_f(r))+O(\max_{1\le i \le q}T_{a_i}(r)).
\end{align*}
Repeating again the argument in the proof of assertion (a), we see that the above inequality holds for all $r\in\R^+$ outside a finite Borel measure set. Then the assertion (b) is proved. \hfill$\square$

\vskip0.2cm
\noindent
{\bf Proof of Theorem \ref{1.2}.}
We denote by $\mathcal I$ the set of all permutations of $q$-tuple $(1,\ldots,q)$. For each element $I=(i_1,\ldots,i_q)\in\mathcal I$, we set
\begin{align*}
N_I&=\{r\in\R^+;N^{[k]}_{(f,a_{i_1})}(r)\le\cdots\le N^{[k]}_{(f,a_{i_q})}(r)\}.
\end{align*}
We now consider an element $I$ of $\mathcal I$, for instance it is $I=(1,...,q)$. Then there is a maximal linearly independent subfamily of the set $\{(f,\tilde a_i);1\le i\le n+1\}$ which is of exactly $k+1$ elements and contains $(f,\tilde a_1)$.  We assume that they are $\{(f,\tilde a_{i_j});1=i_1<\cdots <i_{k+1}\le n+1\}$. For each 
$1\le j\le k+1$, we set $J=\{i_1,...,i_{k+1}\}$
$$ V_j=\left \{i\in\{1,\ldots ,q\}\ ;\ (f,\tilde a_j)\in\biggl ( (f\tilde a_{i_s}); 1\le s\le k+1, s\ne j\biggl)_{\mathcal R}\right \}. $$
Since the space $\biggl ( (f\tilde a_{i_s}); 1\le s\le k+1, s\ne j\biggl)_{\mathcal R}$ is of dimension $k$, the set $V_j$ has at most $n$ elements. Hence 
$$\sharp \bigcup_{j=1}^{k+1}V_j=\sharp \bigcup_{j=1}^{k+1}(V_j\setminus J)+(k+1)\le (n-k)(k+1)+(k+1)=(n-k+1)(k+1).$$
Therefore, there exists an index $i_0\le (n-k+1)(k+1)+1$ such that $i_0\not\in\bigcup_{j=1}^{k+1}V_j$. This yields that the set $\{(f,\tilde a_{i_j}); 0\le j\le k+1\}$ is  minimal over $\mathcal R$. Then by Lemma \ref{3.1}, for all $r\in N_I$ we have
\begin{align*}
||\ T(r,f)&\le\sum_{j=0}^{k+1}N^{[k]}_{(f,a_i)}(r)+o(T_f(r))+O(\max_{1\le i \le q}T_{a_i}(r))\\ 
&\le N^{[k]}_{(f,a_1)}(r)+\sum_{i=n-k+2}^{n+1}N^{[k]}_{(f,a_i)}(r)+N^{[k]}_{(f,a_{(n-k+1)(k+1)+1})}(r)+o(T_f(r))+O(\max_{1\le i \le q}T_{a_i}(r))\\
&\le\dfrac{1}{n-k+1}\left (\sum_{i=1}^{n-k+1}N^{[k]}_{(f,a_i)}(r)+\sum_{i=n-k+2}^{(n-k+1)(k+1)}N^{[k]}_{(f,a_i)}(r)+\sum_{i=(n-k+1)(k+1)+1}^{(n-k+1)(k+2)}N^{[k]}_{(f,a_i)}(r)\right )\\
&+o(T_f(r))+O(\max_{1\le i \le q}T_{a_i}(r))\\ 
&=\dfrac {1}{n-k+1}\sum_{i=1}^{(n-k+1)(k+2)}N^{[k]}_{(f,a_i)}(r)+o(T_f(r))+O(\max_{1\le i \le q}T_{a_i}(r))\\
&\le\dfrac {1}{n-k+1}\cdot\dfrac{(n-k+1)(k+2)}{q}\sum_{i=1}^{q}N^{[k]}_{(f,a_i)}(r)+o(T_f(r))+O(\max_{1\le i \le q}T_{a_i}(r))\\
&=\dfrac{k+2}{q}\sum_{i=1}^{q}N^{[k]}_{(f,a_i)}(r)+o(T_f(r))+O(\max_{1\le i \le q}T_{a_i}(r)).
\end{align*}
Repeating again the argument in the proof of Theorem \ref{1.1}, we see that the above inequality holds for all $r\in\R^+$ outside a finite Borel measure set. Hence, the theorem is proved. \hfill$\square$

\vskip0.2cm
{\footnotesize 
\noindent
{\sc Si Duc Quang}\\
Department of Mathematics, Hanoi National University of Education,\\
136-Xuan Thuy, Cau Giay, Hanoi, Vietnam.\\
\textit{E-mail}: quangsd@hnue.edu.vn

\end{document}